\documentclass[english,twocolumn,superscriptaddress,a4paper]{revtex4}

\usepackage[T1]{fontenc}
\usepackage[latin9]{inputenc}
\usepackage{color}
\usepackage{babel}
\usepackage{array}
\usepackage{textcomp}
\usepackage{amsmath}
\usepackage{graphicx}
\usepackage[unicode=true,
 bookmarks=false,
 breaklinks=false,pdfborder={0 0 1},backref=false,colorlinks=true]
 {hyperref}
\hypersetup{
 hidelinks,hidelinks,citecolor=blue,hidelinks,citecolor=blue}

\makeatletter

\begin{document}

\title{Equilibration of energy in slow-fast systems}

\author{Kushal Shah}
\email{kkshah@ee.iitd.ac.in}
\affiliation{Dept of Electrical Engineering, Indian Institute of Technology (IIT)
Delhi, New Delhi 110016, India.}

\author{Dmitry Turaev}
\email{dturaev@imperial.ac.uk}
\affiliation{Dept of Mathematics, Imperial College, London SW7 2AZ, United Kingdom,
and Lobachevsky University of Nizhny Novgorod, 603950 Russia.}

\author{Vassili Gelfreich}
\email{v.gelfreich@warwick.ac.uk}
\affiliation{Mathematics Institute, University of Warwick, Coventry CV4 7AL, United
Kingdom.}

\author{Vered Rom-Kedar}
\email{vered.rom-kedar@weizmann.ac.il}
\affiliation{The Estrin Family chair of Computer Science and Applied Mathematics, The Weizmann Institute
of Science, Rehovot 76100, Israel.}

\begin{abstract}
Ergodicity is a fundamental requirement for
a dynamical system to reach a state of statistical equilibrium. On
the other hand, it is known that in slow-fast systems
ergodicity of the fast sub-system impedes the equilibration of the whole
system due to the presence of adiabatic invariants. Here, we show that
the violation of ergodicity in the fast dynamics effectively drives the whole system
to equilibrium. To demonstrate this principle we investigate dynamics of the so-called
springy billiards. These consist of a point particle of
a small mass which bounces elastically in a billiard where one
of the walls can move - the wall is of a finite mass and is attached
to a spring. We propose a random process model for the slow wall dynamics
and perform numerical experiments with the springy billiards themselves and the model.
The experiments show that for such systems equilibration is always
achieved; yet, in the adiabatic limit, the system equilibrates with a positive exponential rate
only when the fast particle dynamics has more than one ergodic component for certain
wall positions.
\end{abstract}


\maketitle


\section{Introduction}
In classical statistical mechanics, one deals with systems of a large
number of degrees of freedom, where the exact knowledge of the state
of the system at any given moment of time is impossible or irrelevant.
One, therefore, declares the state variables ``microscopic'' and
tries to describe the statistics of ``macroscopic'' variables
(certain functions of the microscopic state) in an ensemble of many
systems similar to the given one. In other words,
statistical mechanics examines averaging over the phase space
of a dynamical system, typically a Hamiltonian one. There is no a
priori way of choosing the probability distribution over which the
averaging is performed: given a dynamical system, its evolution is
different for different initial conditions, and the distribution of
the initial conditions is not encoded in the system and can be arbitrary.
However, one notices that the microscopic variables are usually changing
fast, so the observed macroscopic quantities are, in fact, time averages.
If the system is ergodic with respect to the Liouville measure (the
uniform measure in the phase space restricted to a given energy level),
then the Birkhoff ergodic theorem allows one to replace time averaging
by averaging with respect to the Liouville measure. In other words,
ergodicity dictates that the averaging must be universally performed
over the micro-canonical ensemble \cite{Huang,Reif,Dorfman}. After this choice
of the ensemble is made, standard results of statistical mechanics
are recovered (e.g. for a system of a large number of weakly coupled
systems, each of which has a bounded energy, the averaging over the
Liouville measure in the phase space of the full system yields the
canonical Gibbs distribution of the energies of the constituent systems \cite{Hil14}).

The main problem is that Hamiltonian systems are usually not ergodic,
even if the number of degrees of freedom is large. For example, the
gas of hard spheres is, most probably, ergodic \cite{Sinai1979,Simanyi1999},
but replacing the collisions of the spheres by mutual repulsion will,
quite probably, ruin the ergodicity, even for an arbitrarily steep repulsing
force potential \cite{Rapoport2006,Kaplan2001}. In general, the picture of dynamics
in a smooth potential appears to be of a chaotic sea (a hyperbolic
set in the phase space) with stability islands (regions in the phase
space that contain a positive measure set filled by KAM-tori) 
\cite{Zas85,lict92}. When the islands occupy a noticeable portion of the phase space, the
use of the micro-canonical ensemble for averaging is unfounded. This
problem can be lifted by postulating (as, for example, an experimental
fact) that the systems of physical interest are sufficiently close
to ergodicity (the islands are small).


However, such universal ``apparent ergodicity'' postulate contains
an intrinsic flaw. Indeed, consider an isolated system with a
few slow degrees of freedom, the rest being fast. Would the fast subsystem
be universally ergodic, the evolution of slow variables would obey
adiabatic laws. In this case the full system would have a conserved
quantity other than the energy - the Gibbs volume entropy of the fast
subsystem as a function of the slow variables (one can view
the slow variables as parameters of the fast system that change adiabatically,
and adiabatic processes are known to keep the entropy constant \cite{Hil14}). One
has to have significant fluctuations in the fast subsystem in order
to destroy this additional conserved quantity, otherwise the full slow-fast
system will not appear ergodic on a long time scale. A rigorous formulation
of this fact is given by Anosov-Kasuga averaging theorem \cite{AnosovKasuga}.
A celebrated example of the ergodicity in the fast subsystem preventing
the equilibration in the full system is the ``notorious piston problem'':
it immediately follows from the adiabatic compression law that the system of two ideal gases at
different temperatures contained in a finite cylinder and separated
by an adiabatic movable piston never comes to equilibrium, which seems to defy
the second law of thermodynamics \cite{LebPiSi00,Lieb99,Neishtadt2004,Wright2007}.

In this paper we resolve this issue by proposing a general mechanism
for the onset of an apparent phase space ergodicity and mixing in
slow-fast Hamiltonian systems. It is not based on the assumption of
a large number of degrees of freedom, nor on the inherent instability
of dispersing geometries (such as the hard spheres models). Instead,
we assume that the fast subsystem is not ergodic for a significant
range of values of the slow variables. We call such systems multi-component,
as the fast subsystem has several ergodic components on its energy
level. For simplicity, we also assume that for some values of the
slow variables the fast subsystem is ergodic (i.e., has only one ergodic
component). We demonstrate that in this case slow observables converge exponentially to the vicinity
of their averages with respect to the
Liouville measure for the full system. This suggests, somewhat paradoxically,  that the non-ergodic
behaviour of the fast degrees of freedom leads to equilibration of the full system.

To elucidate this principle we construct a specific realization of
slow-fast systems - the springy billiard models. This is a point particle
of a small mass $m$ which bounces elastically in a billiard in which
one of the walls, hereafter called the bar, can move. The bar is heavy (has mass $M\gg m)$
and is attached to a spring, so, typically, the bar motion is slow
and the particle motion is fast.

When the mass ratio $m/M$ vanishes, the bar motion is independent
of the particle dynamics and the particle gains or loses speed at
the collisions with the bar without changing the periodic motion of
the bar. This process is studied under the name of Fermi acceleration
\cite{Lichtenberg1980,Jarz1993,Lieberman1998,Loskutov2000,Cohen2003,Leonel2006,Lenz2008,Dolgopyat2008,Itin2012,Shah2010,Shah2011,Gelfreich2011,Gelfreich2012,Gelfreich2014,Batistic2014,Shah2015,Pereira2015,Turaev2016}.
It has been established that there are two types of Fermi acceleration.
If the frozen billiard is ergodic for all possible positions of the
bar (ergodic accelerators), then the ensemble averaged kinetic energy
may grow at most quadratically in time \cite{Loskutov2000,Leonel2006,Gelfreich2011,Gelfreich2012}.
If the ergodicity of the frozen
billiard is broken for a part of the period of the bar motion (a multi-component
accelerator), then the particle kinetic energy growth is exponential
in time for almost every initial condition \cite{Gelfreich2011,Shah2010,Shah2011,Shah2015,Batistic2014,Gelfreich2014,Pereira2015}.
Figure \ref{fig:Billiard} presents the billiard tables we consider
here,  including both ergodic (stadium) and multi-component (rectangle with a bar, mushroom) accelerators.

In the finite mass case ($m/M>0$), the particle energy must remain
bounded for all time, and the question of acceleration is replaced
by the question of equilibration - will the kinetic energies of the
particle and the moving bar equilibrate? The numerics we perform always show exponential
equilibration, however there is a sharp difference between the ergodic
and multi-component cases. The equilibration time within the ergodic
accelerator tends to infinity as $m/M\to0$, while it stays bounded
in the multi-component case.

We analyze this effect by deriving a model for the slow bar motion
where the force exerted on the bar by the particle is found by averaging
over one of the ergodic components. In the multi-component case, the
change in the slow variables (the bar position) drives the fast system
(the particle in the instantaneous billiard) to switch between
different ergodic components. Assuming that the switching is well approximated
by a Markov process, we derive  the switching probabilities
and construct a stochastic model for the bar motion (with no adaptable parameters).
We verify the model by numerically comparing its behavior with the behavior of the
corresponding springy billiards at small $m/M$. As expected, in the
ergodic case the model for the slow bar motion has a conserved quantity,
so the equilibration in the corresponding springy billiard occurs
only due to fluctuations from the averaged motion, which explains
the slowing down of the equilibration as the separation of time scales
increases. In the multi-component case, the Markov
process of hopping between the ergodic components leads to equilibration
at positive rates. Numerics shows that the bar motion is quite accurately
represented by this model and the equilibration rates remain non-zero
and get close to the model Markov process rates as $m/M\to0$.

The derivation of stochastic model is not billiard specific. A similar Markov process can be constructed for an arbitrary slow-fast
system with a non-ergodic fast subsystem, cf. \cite{Pereira2015}. Under standard irreducibility and aperiodicity conditions, such a process
should converge to a unique stationary measure, and by uniqueness this must correspond to the Liouville measure of the full system. Therefore, we believe the
proposed apparent ergodisation mechanism should be universally applicable.

\section{A Particle in a springy billiard}

\begin{figure}
\begin{centering}
\includegraphics[scale=0.3]{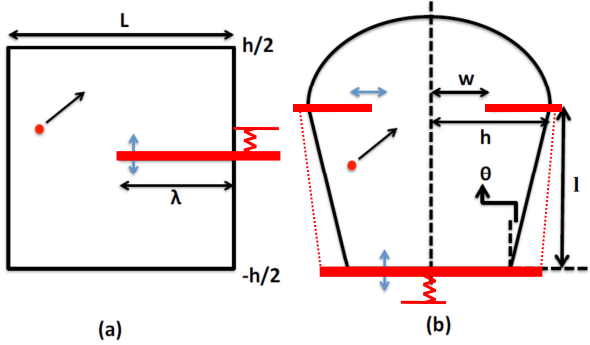}
\par\end{centering}
\caption{Three springy billiards: a particle with a small mass $m$ bounces elastically
within a bounded domain, where one of the walls, the oscillating bar,
has a mass $M=1$ and a spring constant $k$. (a) The  \textbf{Rectangle with
Oscillating Bar} (ROB) model. Numerical values in all simulations are: $\lambda=1,L=h=2,k=81,u_{p}=18/\sqrt{5},E=1,$
so  $y_{b}^{max}=\sqrt{2}/9$
(b) The  \textbf{Slanted
Stadium} ($w=1$) and the \textbf{Slanted Mushroom} $w(y_{b})\leq h$ (Eq. \ref{eq:w_neck}). Numerical values in all simulations are: $E=1,h=1,\tan\theta=0.17$, $\ell=\ell_{0}-y_{b}$, $\ell_{0}=2$,
$k=1$ so $y_{b}^{\max}=\sqrt{2}$.  \label{fig:Billiard}}
\end{figure}

\begin{figure*}
\begin{centering}
\includegraphics[scale=0.85]{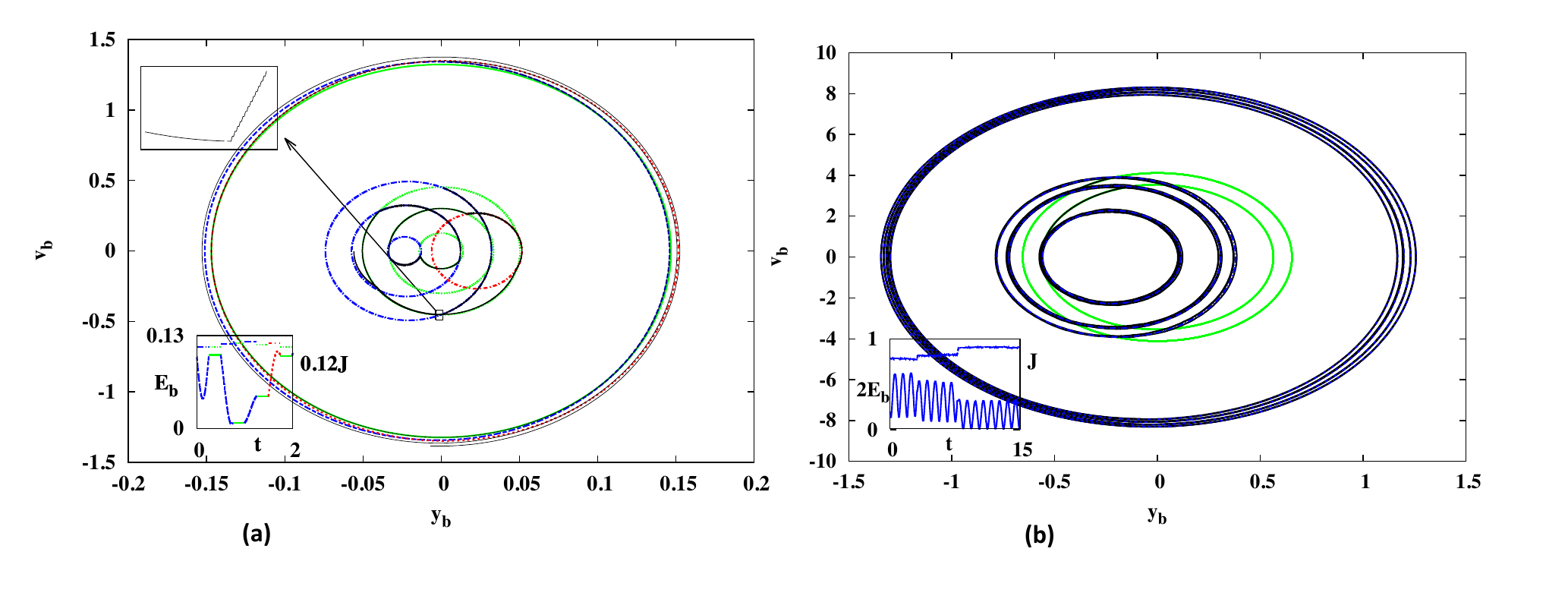}
\par\end{centering}
\caption{Piecewise adiabatic bar motion for the (a) ROB and (b) mushroom springy
billiards. Trajectories of the springy billiard (black lines) for $m=10^{-8}$ with $E_{p}(0)=0.9$ (at the center) and  $E_{p}(0)=0.1$ (at the periphery) follow
closely the corresponding level lines of the piecewise adiabatic invariants (green - level lines of the free bar motion, blue/red
- level lines of the bar motion for a particle hitting it from above/below). The
lower insets show the oscillations and piecewise constant form of $E_{b}(t)$
and $J(t)$ respectively, the upper inset shows the transfer from one level set to another. \label{fig:Level-Sets}}
\end{figure*}








Consider a mass $m\ll1$ particle in a $d$-dimensional billiard. One
of the billiard walls, the bar of mass $M=1$, is suspended on a spring
with a spring constant $k$, so it may oscillate vertically. At impact,
the bar and the particle undergo an elastic collision leading to the
exchange of momentum and energy between the bar and the particle:
\[
v_{p}'=\frac{2v_{b}-\left(1-m\right)v_{p}}{1+m}\,,\:\qquad v_{b}'=\frac{\left(1-m\right)v_{b}+2mv_{p}}{1+m}\,,
\]
where $v_{b}$, $v_{p}$ and $v_{b}'$, $v_{p}'$ are the vertical velocities
of the bar and the particle just before and just after the collision.
The total energy of the system $E=E_{p}+E_{b}$ is preserved, whereas
the particle energy $E_{p}=m\frac{\mathbf{v}_{p}^{2}}{2}$ (hereafter
$\mathbf{v}_{p}$ denotes the particle velocity) and the bar energy
$E_{b}=\frac{v_{b}^{2}+ky_{b}^{2}}{2}$ ($y_{b}$ is the bar position)
change at impact.

The system ``bar-particle'' is a slow-fast $(1+d)$ degrees of freedom
Hamiltonian system: $|v_{b}|\le\sqrt{2E}$ is of order one whereas
the particle speed is typically large ($\|\mathbf{v}_{p}\|=\sqrt{2E_{p}{}/m}$).
Usually, when the particle moves sufficiently fast, many collisions
with the bar occur in short time intervals, so the averaged motion
of the bar is governed by the equation
\begin{equation}
\ddot{y_{b}}+U'\left(y_{b}\right)=F(y_{b},E_{p}),\label{eq:ybdd}
\end{equation}
where $U\left(y_{b}\right)=\frac{1}{2}ky_{b}^{2}$ is the spring potential,
and $F(y_{b},E_{p})$ denotes the averaged force exerted on the bar at position
\(y_{b}\) by the particle with the energy $E_{p}$. The averaging is performed
over many collisions at a frozen value of $y_{b}$. Since the work
done by this force corresponds to the change in the particle's energy,
we conclude that $F=-dE_{p}\big/dy_{b}$.

If the frozen billiard is ergodic for each value of $y_{b}$,
then the Anosov-Kasuga theorem \cite{AnosovKasuga,Jarz1993,Neishtadt2004,Wright2007,Gelfreich2012}
implies that the phase space volume under a given energy level is
approximately preserved (for most trajectories, for sufficiently small
$m$ and on any finite interval of the slow time), hence
\begin{equation}
J=E_{p}^{d/2}V_{c}(y_{b})\approx\mathrm{const},\label{eq:adiabaticlaw}
\end{equation}
where $V_{c}(y_{b})$ is the volume of the billiard domain. This implies
that the average force acting on the bar equals to $F=\frac{2}{d}E_{p}(y_{b})\frac{V_{c}'(y_{b})}{V_{c}(y_{b})}$.
Note that the same formula follows from the ideal gas law. Since $E_{p}=E-E_{b}=E-\frac{1}{2}\dot{y}_{b}^2-U(y_{b})$,
Eq.~\ref{eq:ybdd} becomes
\begin{equation}
\ddot{y_{b}}+U'\left(y_{b}\right)=\frac{2}{d}[E-\frac{1}{2}\dot{y}_{b}^{2}-U(y_{b})]\frac{V_{c}'(y_{b})}{V_{c}(y_{b})},\label{ergbdd}
\end{equation}
One may check that the bar motion defined by Eq. \ref{ergbdd} indeed
follows the level sets of $J$. By noting that the adiabatic law of Eq. \ref{eq:adiabaticlaw}
takes the form $\left(\frac{J}{V_{c}(y_{b})}\right)^{2/d}=E_{p}$,
we find that the adiabatic bar motion is governed by an effective
potential
\begin{equation}
U_{\text{eff}}(y_{b})=U(y_{b})+\left(\frac{J}{V_{c}(y_{b})}\right)^{2/d}\label{efpot}
\end{equation}
where $J$ is determined by the initial condition. Similar effective
potential was derived in the context of the piston problem \cite{Neishtadt2004,Wright2007}.
It is easy to see that $0\le J\le J_{f}=(E-U(y_{f}))^{2/d}V_{c}(y_{f})$.
Here $y_{b}=y_{f}$ corresponds to the pressure equilibrium, where
the pressure due to the collisions with the particle is compensated
by the force exerted by the spring ($y_{f}$ depends on
the total energy $E$). At the other extreme of $J=0$, the particle
does not move and all the energy is in the oscillating bar.

\begin{figure}
\begin{centering}
\includegraphics[scale=0.32]{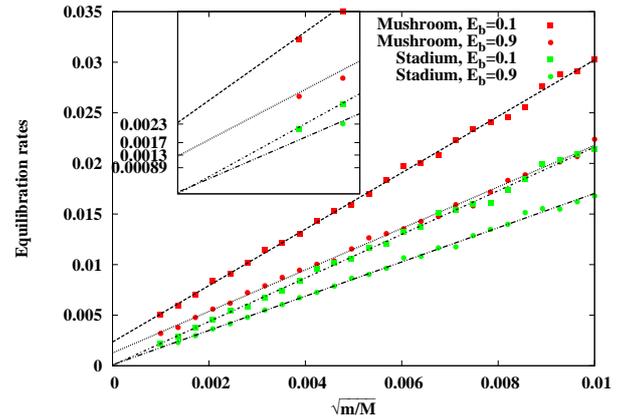}
\par\end{centering}
\caption{Energy equilibration rates dependence on  \(\sqrt{m/M}\) (ensembles of 10000
particles in each computation) for the springy mushroom and stadium.
The linear fit to the springy stadium equilibration rates produces, when extrapolated to $m\big/M=0$, quite small values: $-0.000003$ for $E_b=0.1$ and
$0.000090$ for $E_b=0.9$. The linear extrapolation of springy mushroom equlibration rates to $m\big/M=0$ gives $0.0023$ for $E_b=0.1$ and $0.0013$ for $E_b=0.9$;
these values are comparable to the positive rates computed for the stochastic model,  Eq. \ref{eq:Mushroom-Eqn}.
\label{fig:EquilibrationRate}}
\end{figure}

We performed numerical studies for a classical example of ergodic
billiard - the slanted half-stadium \cite{BuniStad}, with the springy oscillating bar at the bottom, see Fig. \ref{fig:Billiard}b
with $w=1$.
When the initial speed of the particle in the stadium is large enough,
we observe that the bar motion indeed follows a level line of $J$
for many bar oscillations. Since for all possible values of $J$ and
$E$ the motion in the potential given by Eq. \ref{efpot} is periodic, as long
as the adiabatic invariant $J$ stays nearly constant the full system
does not equilibrate. Notably, we demonstrate numerically
that the motion for $m>0$ does lead to equilibration
on a sufficiently long time scale. We also note that in accordance
with the adiabatic theory, the numerically found equilibration time
tends to infinity as $m\to0$ (see the next section and Fig. \ref{fig:EquilibrationRate}).


Next, we present two \textbf{multi-component} springy billiards, in which the ergodicity of the fast dynamics is broken for
some intervals of time. The first is the Rectangle with Oscillating
Bar (ROB) \cite{Shah2010}, Fig. \ref{fig:Billiard}a. A particle moves in a rectangle (of length
$L$ and height $h$) which is partially split by a length
$\lambda$ horizontal bar attached to a spring. The particle horizontal
speed $|u_{p}|$ is preserved, so the horizontal motion is periodic
with period $T_{p}=2L/|u_{p}|$. This period is divided into two time
intervals. On the first one, the particle does not hit the bar and
its vertical speed is preserved. On the second interval, of length
$\tau T_{p}$, where $\tau=\lambda/L$, the particle enters a chamber
above or below the bar where it gains or lose vertical speed as it
hits the moving bar many times. Consequently, during a single passage
above or below the bar, the vertical speed $|v_{p}|$ obeys the
adiabatic law given by Eq. \ref{eq:adiabaticlaw} with $d=1$, $E_{p}=\frac{1}{2}mv_{p}^{2}=E-\frac{1}{2}v_{b}^{2}-\frac{1}{2}ky_{b}^{2}$
and $V_{c}^{up}(y_{b})=1-y_{b}$ when the particle is above the bar
and similarly, when below, $V_{c}^{down}(y_{b})=1+y_{b}$ (here \(E, E_{p}\) do not include the particle horizontal kinetic energy which is decoupled from the dynamics).
Hence $J_{up/down}=\sqrt{E-\frac{1}{2}v_{b}^{2}-\frac{1}{2}ky_{b}^{2}}(1\mp y_{b})=const$. We further assume that the deterministic process is
well approximated by a stochastic one, by which the probabilities
to enter the chamber above/below the bar are proportional to the length
of the gap between the bar and the upper or, respectively, lower boundary
of the rectangle, i.e., these probabilities are equal to $(1\mp y_{b})/2$
where $y_{b}$ is taken at the moment of entrance to the chamber.
The same assumptions were used in \cite{Shah2010} for the study of
Fermi acceleration in ROB in the case of infinitely heavy bar (the
limit $m=0$); as numerics performed in \cite{Shah2010} show, these
assumptions hold in the case $m=0$ with a good precision.

\begin{figure}
\begin{centering}
\includegraphics[scale=0.7]{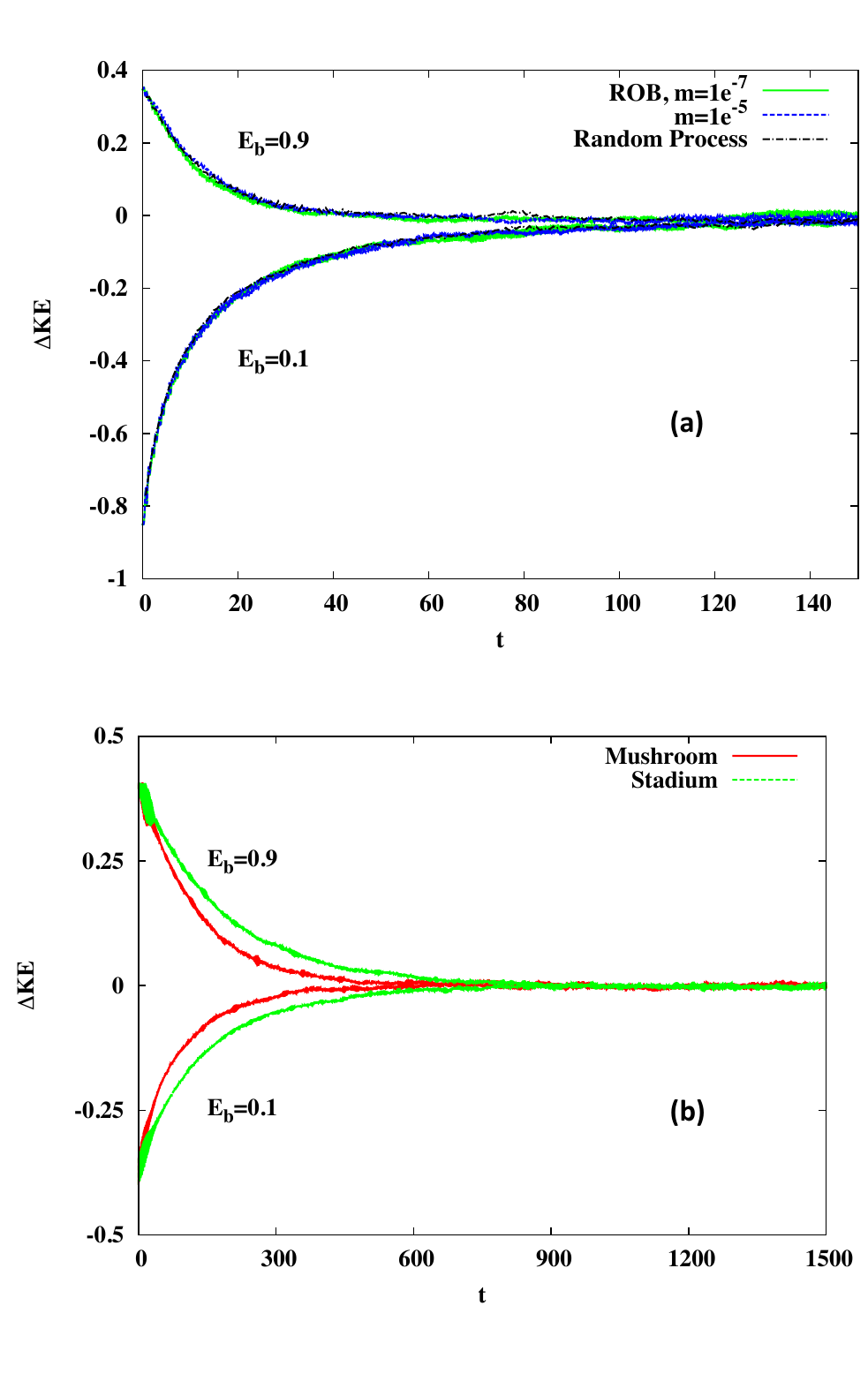}
\par\end{centering}
\caption{Kinetic energy equilibration  for springy billiards. (a) 
The ROB equilibration depends on the initial ensemble and is  independent of $m$: the equilibration processes for ensembles of 6000 particles starting at $E_{b}(0)=0.9$
(and similarly for $E_{b}(0)=0.1$) are essentially the same for two mass ratios ($m/M=10^{-5}$ and $10^{-7}$) and for the random process model, Eq. \ref{eq:Wall-Eqn-Up}.
(b) The stadium and mushroom equilibration process for $m/M=10^{-5}$. The equilibration rates depend on both the initial ensemble and the mass ratio, see Fig. 4. In all cases
the stadium equilibrates slower than the mushroom.
\label{fig:Equilibrationproc}}
\end{figure}

Since ${E_{p}}=J_{up/down}^{2}/(1\mp y_{b})^{2}$, the right-hand side of Eq.~\ref{eq:ybdd} becomes $\mp2J_{up/down}^{2}/(1\mp y_{b})^{3}$
for particles above/below the bar, respectively. Hence, we suggest
that the bar-particle system is well approximated by the following
$T_{p}$-periodic probabilistic hybrid system:
\begin{equation}
\ddot{y_{b}}+ky_{b}=\!\left\{\!\!\!\!\!\! \begin{array}{cc}
\begin{array}{clc}
-\frac{2J_{j,up}^{2}}{(1-y_{b})^{3}} &\! \text{prob.} &\!\!\! \beta_{j}\\
\frac{2J_{j,down}^{2}}{(1+y_{b})^{3}} &\! \text{prob. } &\!\!\! 1-\beta_{j}
\end{array}
& \left\{ \frac{t}{T_{p}}\right\} \in[0,\tau)\\
0 & \left\{ \frac{t}{T_{p}}\right\} \in[\tau,1),
\end{array}\right.\label{eq:Wall-Eqn-Up}
\end{equation}
$j=\left\lfloor \frac{t}{T_{p}}\right\rfloor $, $J_{j,up/down}=\sqrt{E-E_{b}(jT_{p})}\left(1\mp y_{b}(jT_{p})\right)$,
$E_{b}=\left(\dot{y}_{b}^{2}+ky_{b}^{2}\right)/2$ and $\beta_{j}=(1-y_{b}(jT_{p}))/2$.
Namely, the bar dynamics follows the level lines of $J_{up}$, $J_{down}$,
or $J_{0}=\frac{1}{2}(v_{b}^{2}+ky_{b}^{2})$. At $t=jT_{p}$, the
bar height $y_{b}(jT_{p})$ determines the probability $\beta_{j}$
to choose $J_{up}$ vs. $J_{down}$, and $(y_{b}(jT_{p}),\dot{y}_{b}(jT_{p}))$
determine the particular $J_{up/down}$ level set along which the
motion will continue. Fig. \ref{fig:Level-Sets}a demonstrates that
the bar motion of the springy ROB follows the level sets of the corresponding
$J$'s. In the next section, we examine equilibration, demonstrating
that at small $m$ the stochastic model given by Eq. \ref{eq:Wall-Eqn-Up} provides
a good approximation to bar dynamics in the springy ROB billiard.

The ROB model represents a family of multicomponent billiards in which
the transition from a single ergodic component to two (or more) possible
ergodic components and back to a single component occurs at discrete
sequence of prescribed times (e.g. at $\{t/T_{p}\}\in\{0,\tau\}$
in the ROB model). Similar behavior occurs when one splits and unites
$d$-dimensional chaotic billiards by a moving wall (see \cite{Gelfreich2011,Gelfreich2012}).
In particular, for such systems the standard adiabatic law applies
at each interval of time along which the particle motion corresponds
to a single ergodic component, and a stochastic process similar to
Eq. \ref{eq:Wall-Eqn-Up} provides a good model for the bar dynamics.

Next, we describe a different multi-component system,  the slanted springy mushroom of Fig. \ref{fig:Billiard}b \cite{Gelfreich2014, BuniMush}, in which particles
leak from one ergodic component to another during some time interval.
This is expected to be the generic behavior in systems with mixed
phase space (see also \cite{Batistic2014,Pereira2015}). The base of the mushroom stem is a springy bar with vertical position $y_{b}$, so the mushroom
stem length is $\ell=\ell_{0}-y_{b}$. The mushroom throat width
$w$, which controls the capture and release of particles from the
mushroom cap, changes in synchrony with the bar position as prescribed
later on. This is a multi-component system - the motion is mixing
in the stem and integrable in the cap \cite{BuniMush}.
The oscillations in $w$ lead to the exchange of particles between
the integrable and chaotic zones \cite{Gelfreich2014}. In such situation,
by which a particle in the chaotic component may leak into a different
component at any instant of time in some time interval, the usual
adiabatic law needs to be modified to the ``leaky adiabatic law'':
\begin{equation}
\frac{dE_{p}}{E_{p}}=-\frac{dV}{V_{c}}\label{eq:leakyadiab}
\end{equation}
(see \cite{Gelfreich2014} for derivation and corroboration for $m=0$).
Here, $V$ is the total phase space volume for the frozen mushroom
billiard on the energy level $E_{p}=1$ and $V_{c}$ is the phase
space volume of the chaotic zone on the same energy level. Then $V_{c}=V-V_{ell}$
(see details in \cite{Gelfreich2014}) where 
$V(y_{b})=\pi^{2}h^{2}+2\pi\left(2h\ell-\ell^{2}\tan\theta\right),
V_{ell}(y_{b})=V_{ell}(w(y_{b}))
=2\pi h^{2}\left(\cos^{-1}\frac{w}{h}-\frac{w}{h}\sqrt{1-\left(\frac{w}{h}\right)^{2}}\right).$ 
This leaky adiabatic law implies that as long as the particle remains in
the chaotic component the system has an adiabatic invariant given
by $J=E_{p}(y_{b},\dot{y}_{b})G(y_{b})$ where $G(y_{b})=\exp\left(\int_{0}^{y_{b}}\frac{V'(s)}{V_{c}(s)}ds\right).$
Hence, while the particle is not captured in the cap, the bar motion
occurs along the level lines of $J=\left(E-\frac{1}{2}\dot{y}_{b}^{2}-\frac{1}{2}ky_{b}^{2}\right)G\left(y_{b}\right),$
producing an effective potential: $U_{\text{eff}}(y_{b})=U(y_{b})+\frac{J_{0}}{G(y_{b})}$
where $J_{0}=J(y_{b}(0),\dot{y}_{b}(0))$. When the elliptic zone
is expanding $(\dot{w}<0)$, the particle may transfer from the chaotic
to the elliptic zone. Since the motion in the stem is chaotic, we
model the capture time as a random variable. We set the probability
to transfer at the time interval $\left(t_{c},t_{c}+dt\right)$ from
the stem to the cap as the ratio of the transferred phase space volume
to the chaotic zone volume (see \cite{Gelfreich2014}):
\begin{eqnarray}
P_{cha\rightarrow ell}(t_{c})|_{\left\{ t_{c}|\dot{w}<0\right\} }=\frac{V_{ell}'(y_{b})\dot{y}_{b}}{V_{c}\left(y_{b}\right)}dt=d\ln\frac{G(y_{b})}{V_{c}(y_{b})}.\label{eq:P_cha-ell}
\end{eqnarray}
Once captured in the cap, the particle does not influence the bar,
so the bar and particle energies are preserved separately and the
bar moves only due to the spring force. The particle gets released
from the cap at the time interval $\left(t_{r},t_{r}+dt\right)$ where
$t_{r}$ is the first instance at which the throat reached back the
width it had at the time of capture: $w(t_{r})=w(t_{c})$ and in the
adiabatic limit $dt$ vanishes. If $\ell$ is not a single-valued
function of $w$, the bar position is changed between the capture
and the release time: $y_{b}(t_{c})\neq y_{b}(t_{r})$, so after the
capture episode the particle follows a new level set of $J=J_{r}=J(y_{b}(t_{r}),\dot{y}_{b}(t_{r}))$.
To ensure that $\ell$ is not single valued function of $w$ for all
oscillations of the bar, we choose a protocol for the dependence of
$w$ on $\ell$ for which $w$ takes its minimum at $\ell_{f}=\ell_{0}-y_{f}$
(recall that $y_{f}$ is the pressure equilibrium point where $J$
takes its maximal value, here $y_{f}\approx -0.2436$):
\begin{align}
w(y_{b}) & =\begin{cases}
\min\left\{ 1,\:0.7+0.6(y_{b}-y_{f})^{2}\right\}  & \mbox{for }y_{b}<y_{f},\\
\min\left\{ 1,\:0.7+6(y_{b}-y_{f})^{2}\right\}  & \mbox{for }y_{b}\ge y_{f}.
\end{cases}\label{eq:w_neck}
\end{align}
Hence, the stochastic model for the bar motion becomes
\begin{equation}
\ddot{y_{b}}+U'\left(y_{b}\right)=\!\left\{\!\! \begin{array}{ll}
\frac{V'(y_{b})}{V_{c}(y_{b})}\left(1-\frac{1}{2}\dot{y}_{b}^{2}-\frac{1}{2}ky_{b}^{2}\right) & \!\text{(particle in }V_{c})\\
0 & \!\text{(particle in }V_{ell})
\end{array}\right.\label{eq:Mushroom-Eqn}
\end{equation}
where the particle moves from the stem to the cap at time $t_{c}$
with probability $P_{cha\rightarrow ell}(t_{c})$ given by Eq. \ref{eq:P_cha-ell},
and subsequently gets released at $t_{r}=t_{r}(t_{c})$ (defined by
$w(y_b(t_r)) = w(y_b(t_c))$). Figure \ref{fig:Level-Sets}b
demonstrates that the bar-particle simulations adhere to Eq. \ref{eq:Mushroom-Eqn}.
We also verified that the distribution of capture events in the cap
as a function of $y_{b}$ is well approximated by Eq. \ref{eq:P_cha-ell}
and that the release time $t_{r}$ is well approximated by the equal throat-width
rule. In the next section, we show that this stochastic model
provides a good approximation to the equilibration process in the
limit $m\rightarrow0$.

\section{Ensemble equilibration rates }

So far, we have checked that the adiabatic theory, by which one averages
the fast particle motion, enables to predict the evolution of the
averaged bar motion for a finite time in each of the ergodic components
of the fast system. In the case where the fast subsystem is ergodic,
the predicted adiabatic bar motion is periodic, so no equilibration
occurs as long as the adiabatic invariant is preserved. In the multi-component
cases the bar motion follows effectively a random dynamical system
by which the bar motion switches between different laws of motion,
leading to chaotization of orbits, and hence equilibration.

To test equilibration, we examine the behavior of ensemble average of kinetic energies. By Equipatition Theorem,
at equilibrium, the kinetic energy for each degree of freedom should be the same \cite{Hil14}.
We thus define \(\Delta KE=\left\langle \frac{M}{2}v_{b}^{2}\right\rangle -\left\langle \frac{m}{2}v_{p}^{2}\right\rangle \)
for the ROB and  $\Delta KE=\left\langle \frac{M}{2}v_{b}^{2}\right\rangle -\left\langle {E_p}\right\rangle/2 $ for the slanted stadium and mushroom.

Numerically, we find that for finite $m$, in all cases,
for all initial ensembles, \(\Delta KE_{t\to+\infty}\rightarrow0\),  see Fig. \ref{fig:Equilibrationproc}.
This suggests that for finite \(m\) the system gets close to equilibrium in finite time for both ergodic and multi-component springy billiards.

Next, we examine the equilibration rates. We find that the convergence to the equilibrium is exponential, yet the rate of convergence depends on the initial ensemble
and, also, on the choice of the interval of fitting. Hereafter, in order to make a comparison between the rates of convergence in different models possible,
we fix a practical definition of the equilibration rate as the best fitted slope to  $\log \left|\Delta KE\right|$ vs. $t$ on the time interval
\([0,T]\) wher $T$ is defined by \(\log |\Delta KE(0)/\Delta KE(T)|\approx 1\).  For each fixed ensemble of initial conditions we examine how such defined
rate depend on the mass ratio. For the ROB case the rates do not display any significant dependence on $m$ (see Fig. \ref{fig:Equilibrationproc}a).
For the mushroom and stadium cases, the equilibration rates increase
from their limit values proportionally to $\sqrt{m}$. By extrapolation to $m=0$,
we find that the limit equilibration rates for the stadium vanish (Fig. \ref{fig:EquilibrationRate}a)
whereas for the mushroom they are positive (Fig. \ref{fig:EquilibrationRate}b).

This is our main finding, as it reveals the profound difference between
the ergodic and multi-components cases. The vanishing of the equilibration
rates in the ergodic case is in agreement with adiabatic theory. The
positive rates achieved in the multi-component case reflect the efficient mixing induced by the hoping between the ergodic components.

The extrapolation of the data obtained by the simulations at non-zero $m$ can be sensitive to noise and the chice of the fitting procedure.
Therefore, to claim that the $m=0$ limit equilibration rates are strictly positive in the multi-component cases, we need to benchmark them
against their theoretical values. To this aim, we also simulate the stochastic models given by Eqs. \ref{eq:Wall-Eqn-Up} and \ref{eq:Mushroom-Eqn}, and observe that
these models also equilibrate in a similar fashion, with the ensemble-dependent exponential rates close to those
obtained by the springy billiard simulations at small $m$ (Fig. \ref{fig:EquilibrationRate}). 
Thus, we compare the ROB rates for initial ensembles with
$E_{b}=0.9E$ and $E_{b}=0.1E$, computed for 10 runs of 6000 particles each. The stochastic simulations of Eq. \ref{eq:Wall-Eqn-Up} produce the rates
$0.0296\pm0.0026$ for the ensembles with the initial $E_{b}=0.9E$ and $0.0960\pm0.0069$ for the ensembles with the initial $E_{b}=0.1E$,
whereas the actual springy billiard simulations produce the rates, respectively, $0.0307\pm0.0039$ and $0.0961\pm0.0057$ for $m=10^{-7}$ and $0.0281\pm0.0025$
and $0.0862\pm0.0076$ for $m=10^{-5}$.
Similarly, the rates for the stochastic simulations of the mushroom model given by Eq. \ref{eq:Mushroom-Eqn}, computed for 10 runs of 10000 particles each, are
$0.00089 \pm 0.00002,\text{ } 0.00170 \pm 0.00007$
whereas the limit values from Fig. \ref{fig:EquilibrationRate}b
are $0.0013$ and $0.0023$ for $E_{b}=0.9E$ and $E_b=0.1E$, respectively.

We conclude that the stochastic models, with no adaptable parameters,
provide a reasonable approximation to the dynamics observed in our numerical experiments.
These results support our claim that the behavior of slow-fast systems
with a multi-component fast subsystem can be modeled by the random process
of Markov switching between several different equations for slow evolution
obtained by averaging over different ergodic components in the fast phase space.

\section{Discussion}


Springy billiards demonstrate an important principle in slow-fast
Hamiltonian systems: ergodicity of the fast subsystem impedes the
ergodisation in the full slow-fast system whereas its violation can
lead to equilibration. Indeed, we computed the equilibration rates
for ergodic and multi-component springy billiards and demonstrated,
by extrapolation to $m=0$, that the equilibration rates
for the ergodic cases vanish in the limit $m=0$ whereas for the multi-component ones
they are positive (Fig. \ref{fig:EquilibrationRate}). We showed that
the limit behavior of the multi-component case is well approximated
by a random process in which the slow variables follow the averaged
dynamics on each ergodic component of the fast system, and switch
randomly between these different averaged systems, leading to complete
chaotization. We stress that in this random model of multi-component
system both the magnitude and probabilities of the jumps remain finite
as $m\rightarrow0$, hence the chaotization time remains bounded
in this limit.

The observed sensitivity of the equilibration rates
to the initial ensemble energy may be related
to the two non-trivial regimes at $E_{p}=0$ (the Fermi acceleration
limit) and $E_{p}=E_{f}$ (the pressure equilibrium). In both limits
there is no exchange of energy between the slow and fast systems.
Spectral properties of the operator governing the evolution of the
densities are affected by these states.

The Fermi acceleration limit may be studied directly by taking an ensemble
of fast particles with vanishing kinetic energy (e.g. $|v_{p}|=O(m^{-a}),a<1/2,E_{p}(0)=O(m^{1-2a})$),
so that initially the particles hardly influence the bar motion. Such particles accelerate exponentially fast
in multi-component accelerators
whereas in ergodic accelerators they accelerate as a power-law in
time. This leads to distinct dependence of the transient time 
 on $m$: for a multi-component springy billiard the transients are
of order $O(\ln(m))$ whereas for an ergodic springy billiard they are much
longer, of order $O(1/\sqrt{m})$.

Slower particles may resonate and possibly freeze, and this singular behavior may have non-trivial effects on
the equilibration \cite{EckY10}. Such a behavior breaks the slow-fast structure and thus it cannot be treated by the stochastic model
we propose. However, since this phase space region has small volume, its influence on the statistical behavior seems to be negligible.

We also propose that springy billiards may be used to study additional dynamical
phenomena in slow-fast systems. First, in this paper we assumed for simplicity
that the fast system is ergodic for a range of the slow variables
and that this range is realized in every cycle of the slow system.
In particular, in our setup, the pressure equilibrium points were always
destabilized (in the springy mushroom case this motivated the choice
of $w(y_{b})$, see Eq. \ref{eq:w_neck}). Springy multi-component
billiards for which this property is violated are easy to construct,
and thus proper conditions under which such systems still lead to
ergodization need to be formulated and studied. Such cases were considered
in the Fermi acceleration limit for billiards \cite{Batistic2014}
and for smooth homogeneous Hamiltonians \cite{Pereira2015}. Second, we considered a
one degree of freedom slow system with a single
fixed point which is always stable (a single bar moving vertically).
Thus, the dynamics in each of the averaged slow systems is trivially
integrable and periodic. More complicated situations, possibly with
several degrees of freedom may be studied. Third, extended simulations
and analysis of the random models may shed light on the role of islands
and the effect of the singular regimes (the Fermi limit and the pressure
equilibria). Fourth, incorporating the finite $m$ fluctuations into
the random models is a challenging problem, for both the ergodic and
multi-component cases.

The multi-component equilibration mechanism in springy billiards is
expected to appear also in smooth slow-fast multi-component
systems. In fact, the new chaotization mechanism we propose is reminiscent
of the phenomena of adiabatic chaos - chaotization of smooth slow-fast
systems in which the fast dynamics is integrable, yet the structure of the fast
phase space changes as the slow variables are changed \cite{Tennyson1986,Vain1998,Neish1991}.
The new suggestion here is that
this mechanism is universal and is not restricted to fast subsystems
that are integrable. In fact, if some of the ergodic components of
the fast system are chaotic, and, moreover, if there exists a range
of slow variables for which the chaotic ergodic component occupies
a large portion of the fast subsystem phase space the equilibration
may be particularly fast.

Finally, we propose a broader viewpoint on this work. Here, our slow
and fast systems were just two mechanical components and the time-scale
separation stemmed from the mass ratio. In the broader statistical
mechanics context, the fast system governs the motion of many particles
(and is thus high-dimensional) and the slow macroscopic variables
are defined as certain averages over the fast, microscopic system. When
the structure of the fast system changes, for example, from
a gas to a liquid state, one declares that a phase transition occurs.
Usually, for macroscopic values near phase transition, the microscopic
phase space structure is complex, with long-lasting structures in
which the two states coexist. We may say that in this range of the
macroscopic variables the fast system is inherently multi-component.
Therefore, we conjecture that phase transitions may play a central
role in the equilibration process between microscopic and macroscopic
variables (e.g. consider the analog of the notorious piston problem
in a multi-phase gas).

\section*{Acknowledgement}

RK acknowledges support by the Israel Science Foundation (grant 1208/16). KS and DT would like to acknowledge the financial
support and hospitality of Weizmann Institute of Science where a part of this work was done. KS also thanks the SERB-DST, Government of India (File No. : SR/FTP/PS-108/2012)
for financial support. DT acknowledges the support of RSF grant 14-41-00044 to this research, and also thanks the Royal Society and EPSRC.
VG's research was supported by EPRC (grant EP/J003948/1).

\end{document}